\documentclass{PNAStwob}

\usepackage[pdftex]{graphicx}

\usepackage{amssymb,amsfonts,amsmath,phonetic}
\usepackage{color}
\usepackage{ulem}
\usepackage{graphicx}
\usepackage{type1cm}
\usepackage{eso-pic}
\usepackage{color}

\def \emph{ \normalem }

\newcommand{\cS}{\mathcal{S}}
\newcommand{\ph}{\varphi}

\newcommand{\G}{\Gamma}

\newcommand{\bD}{\mathbf{D}}
\newcommand{\Bf}{\mathbf{f}}
\newcommand{\bd}{\mathbf{d}}
\newcommand{\bX}{\mathbf{X}}
\newcommand{\bY}{\mathbf{Y}}
\newcommand{\cM}{\mathcal{M}}
\newcommand{\cA}{\mathcal{A}}

\newcommand{\R}{\mathbb{R}}

\newcommand{\abs}[1]{\left\vert#1\right\vert}

\newcommand{\set}[1]{\left\{#1\right\}}

\newcommand{\wt}[1]{\widetilde{#1}} 
\def \ii{\mbox{\rm \footnotesize{i}}}
\definecolor{darkred}{rgb}{0.6,0,0}
\newcommand{\yl}[1]{{\color{darkred}#1}}     

\contributor{Submitted to Proceedings of the National Academy of
Sciences of the United States of America}
\url{www.pnas.org/cgi/doi/10.1073/pnas.0709640104}
\copyrightyear{2008} \issuedate{Issue Date} \volume{Volume}
\issuenumber{Issue Number}

\begin{document}



\title{Algorithms to automatically quantify the geometric similarity of anatomical surfaces}

\author{
D.M.Boyer\affil{1}{Anthropology \& Archaeology Dept., Brooklyn
College, City Univ.  New York, Brooklyn  NY
11210}\affil{2}{Interdept. Doct. Progr. Anthropol. Sciences, Stony
Brook Univ., Stony Brook  NY 11794} Y.Lipman\affil{4}{Dept. of Computer Science and Applied Mathematics, Weizmann Inst. of Science, Israel}\affil{3}{Computer
Science Dept. and Progr. Appl. \& Comp. Math., Princeton Univ.,
Princeton  NJ 08544} E.$\,$St.Clair\affil{2}{} J.Puente\affil{3}{}
B.A.Patel\affil{2}{} \affil{5}{Anatomical Sciences Dept., Stony
Brook Univ. } T.A.Funkhouser\affil{3}{}
J.Jernvall\affil{2}{}\affil{6}{Ecology \& Evolution Dept., Stony
Brook Univ.} \affil{7}{Developmental Biology Program, Inst. for
Biotechnology, University of Helsinki, PO Box 56, FIN-00014
Helsinki, Finland} \and I.Daubechies\affil{3}{}
\thanks{To whom correspondence should be addressed. E-mail:
ingrid\@math.princeton.edu} }

\contributor{Submitted to Proceedings of the National Academy of
Sciences of the United States of America}

\maketitle

\begin{article}

\begin{abstract} We describe
new approaches for distances between pairs of 2-dimensional surfaces
(embedded in 3-dimensional space) that use local structures and
global information contained in inter-structure geometric
relationships. We present algorithms
to automatically determine these distances as well as geometric correspondences.\\
This is motivated by the aspiration of students of natural science
to understand the continuity of form that unites the diversity of
life. At present, scientists using physical traits to study
evolutionary relationships among living and extinct animals analyze
data extracted from carefully defined anatomical correspondence
points (landmarks). Identifying and recording these landmarks is time consuming and
can be done accurately only by trained morphologists. This renders
these studies inaccessible to non-morphologists, and causes phenomics to lag behind
genomics in elucidating evolutionary patterns. \\
Unlike other algorithms presented for
morphological correspondences our approach does not require any
preliminary marking of special features or landmarks by the user.
It also differs from other seminal
work in computational geometry in that our algorithms are polynomial
in nature and thus faster, making pairwise comparisons feasible for
significantly larger numbers of digitized surfaces.
We illustrate
our approach using three datasets representing teeth and different
bones of primates and humans, and show that it leads to highly
accurate results.
\end{abstract}


\dropcap{T}o document and understand physical and biological
phenomena (e.g., geological sedimentation, chemical reactions,
ontogenetic development, speciation, evolutionary adaptation, etc.),
it is important to quantify the similarity or dissimilarity of
objects affected or produced by the phenomena under study. The grain
size or elasticity of rocks, geographic distances between
populations, or hormone levels and body masses of individuals --
these can be readily measured, and the resulting numerical values
can be used to compute similarities/distances that help build
understanding. Other properties like genetic makeup or gross
anatomical structure can not be quantified by a single number;
determining how to measure and compare these is more involved
\cite{Sarg_et_al2007,Schwartz_Tat1985,MacLeod1999,Bookstein2007}.
Representing the structure of a gene (through sequencing) or
quantification of an anatomical structure (through the digitization
of its surface geometry) leads to more complex numerical
representations; even though these are not measurements allowing
direct comparison with their counterparts for other genes or
anatomical structures, they represent an essential initial step for
such quantitative comparisons. The 1-dimensional, sequential
arrangement of genomes and the discrete variation (four nucleotide
base types) for each of thousands of available correspondence points
help reduce the computational complexity of determining the most
likely alignment between genomes; alignment procedures are now
increasingly automated \cite{Liu_et_al2009}. The resulting, rapidly
generated and massive data sets, analyzed with increasing
sophistication and flexible in-depth exploration due to advances in
computing technology, have lead to spectacular progress. For
instance phylogenetics has begun to unravel mysteries of large scale
evolutionary relationships experienced as extraordinarily difficult by morphologists
\cite{Murphy_et_al2001}.

Analyses of massive developmental and genetic data sets outpace
those on morphological data. The comparative study of gross
anatomical structures has lagged behind mainly because it is harder
to determine corresponding parts on different samples, a
prerequisite for measurement. The difficulty stems from the higher
dimension (2 for surfaces vs. 1 for genomes), the continuous rather
than discrete nature of anatomical objects,$\!$\begin{footnote}
{These differences may seem innocuous but they lead to an
exponential increase in the size of the ``search spaces'' to be
explored by comparison algorithms.}
\end{footnote} and from large shape variations.

In standard morphologists' practice, correspondences are first
assessed visually; then, some (10 to 100, at most) feature points can
usually be defined as equivalent and/or identified as landmarks.
Just as comparisons of tens of thousands of nucleotide base
positions are used to determine similarity among genomes, the
coordinates of these dozens of feature points (or measurements they
define) are used to evaluate patterns of shape variation and
similarity/difference \cite{Mitteroeker_Gunz2009}. However, as
stated in 1936 by G. G. Simpson, the paleontologist chaperon of the
Modern Synthesis in the study of evolution, the ``difficulty in
acquiring personal knowledge'' (\cite{Simpson1936} p. 3) of
morphological evidence limits our understanding of the evolutionary
significance of morphological diversity; this remains true today.
New techniques for generating and analyzing digital representations
have led to major advances (see, e.g.
\cite{Wiley_et_al2005,PollyMacleod2008,Zelditch_et_al2004}), but
they typically still require determinations of anatomical landmarks
by observers whose skill of identifying anatomical correspondences
takes many years of training.

Several groups have sought to determine automatic correspondence
among morphological structures. Existing successful methods
typically introduce an effective dimensional reduction, using, e.g.,
2-D outlines and/or images \cite{Lohmann1983}, or, in one of the few
studies attempting automatic biological correspondences in 3-D as a
method for evolutionary morphologists, ``automatically-detected
crest lines'' \cite{ThirionGourdon1996} on surfaces obtained by
CT-scans to register modern human skulls to each other
\cite{Subsol_et_al1998} or to pre-Neanderthal, {\it Homo
heidelbergensis} skulls \cite{Subsol_et_al2002}; another example is
Wiley et al \cite{Wiley_et_al2005}. Studies using 2-D outlines or
images sacrifice a lot of the original geometric information
available in the 3-D objects on which they are based; such
specifically limited representations cannot easily be incorporated
into other studies. More generally and most importantly, none of
these methods are independent of user input. When outlines or
standard 2-D views are used, precise observations of the 3-D
anatomical structures are required by the trained technician who
creates the outlines or 2-D views \cite{Gaston_Oneill2004,
WEEKS_ET_AL1999}. Several methods for 3-D alignment
use
Iterative Closest Point (ICP) algorithms \cite{Besl1992} that
require observer input to fix an initial guess (then further
improved via local optimization); ICP-generated correspondences can
also have large distortions and discontinuities of shape. In \cite{Memoli04,BBK06}
surfaces are matched by using the Gromov-Hausdorff distances between them,
and applications to several shape analysis problems are given. However,
Gromov-Hausdorff distances are hard to compute, and have to be approximated; the gradient
descent optimization used in practice does
not guarantee convergence to a global (rather than local) minimum.

Determination of correspondences or similarities among 3-D
digitizations of general anatomic surfaces
that is both 1) fully automated
and 2) computationally fast (to handle the large data sets that are
becoming increasingly available as imaging technologies become more
widespread and efficient \cite{Schmidt2010}) is still elusive. Our
aim here is to remedy this by fully automating the determination
of correspondences among gross anatomical structures. Success in
this pursuit will help bring to phenomic studies the rate,
objectivity and exhaustiveness of genomic studies. Large scale
initiatives to phenotype model species after systematically knocking
out each gene \cite{Abbott2010}, as well as analysis of
computational simulations of organogenesis
\cite{Salazar-Ciudad_et_al2010} stand to greatly benefit from
automating the determination of correspondence among, and
measurement of, morphological structures.

In this paper, we describe several new distances between surfaces
that can be used for such fully automated anatomic correspondences,
and we test their relevance for biologically meaningful tasks on
several anatomical dataset examples (high
resolution digitizations of bones and teeth).

The paper is organized as follows. Section 1 gives the mathematical
background for our algorithms: conformal geometry and optimal mass
transportation (also known as Earth Mover's distance). In section 2,
we use these ingredients to define new distances or measures of
dissimilarity, including a generalization to surfaces of the
Procrustes distance. Section 3 presents the results obtained by our
algorithms for three different morphological data sets, and an
application.

No technical advance stands on its own; this paper is no
exception. Conformal geometry is a powerful mathematical tool
(permitting the reduction of the study of surfaces embedded in 3-D
space to 2-D problems) that has been useful in many
computational problems; \cite{Gu_Yau_book} provides an introduction
to both theory and algorithms, with many applications, including the
use of conformal images of anatomical structures, combined with user
prescribed landmarks and/or special features, for registration
purposes, seeking ``optimal'' correspondence  between pairs of
surfaces \cite{Gu2004_cortex}. Earth Mover's distances
\cite{rubner2000} and continuous optimal mass transportation
\cite{haker2004} have been used in image registration and for more
general image analysis and parameterization \cite{dominiz10}; in
\cite{Memoli07}
(quadratic) mass transportation is used to
relax the notion of Gromov-Hausdorff distance. Procrustes distances
for discrete point sets are familiar to morphologists and other
researchers working on shape analysis
\cite{Zelditch_et_al2004,Mitteroeker_Gunz2009}.
The mathematical and
algorithmic contribution of our work is the combination in which
we use and generalize these ingredients to construct novel distance
metrics, paired with efficient, fully automatic algorithms not
requiring user guidance. They open the door to new applications
requiring a large number of distance computations.\vspace{-0.5cm}
\section{1. The mathematical components}
\subsection{Conformal geometry}
$\,$A mapping $\ph$ from one 2-dimensional (smooth) surface $\cS$ to
another, $\cS'$, defines for every point $p\in\cS$ a corresponding
point $\ph(p)\in \cS'$. If the mapping is smooth itself, it maps a
smooth curve $\G$ on $\cS$ to a corresponding smooth curve $\G'$ on
$\cS'$ that is called the \emph{image} of $\G$. Two curves $\G_1$
and $\G_2$ on $\cS$ that intersect in a point $s$ are mapped to
curves $\G_1'$, $\G_2'$ that intersect as well, in $s'=\ph(s)$.
Consider the two (straight) lines $\ell_1$ and $\ell_2$ that are
tangent to the curves $\G_1$ and $\G_2$ at their intersection point
$s$; the \emph{angle} between $\G_1$ and $\G_2$ at $s$ is then taken
to mean the angle between the two lines $\ell_1$ and $\ell_2$;
similarly, the angle between the curves $\G_1'$ and $\G_2'$ (at
$s'=\ph(s)$) is the angle between their tangent lines at $s'$. The
mapping $\ph$ is called \emph{conformal} if \emph{for any two smooth
curves $\G_1$ and $\G_2$ on $\cS$, the angle between their images
$\G_1'$ and $\G_2'$ is the same as that between $\G_1$ and $\G_2$ at
the corresponding intersection point}.

Riemann's {\em uniformization theorem} \cite{Gu_Yau_book} guarantees
that every (reasonable) 2-d surface $\cS$ in our standard 3-d space that is a \emph{disk-type surface}
(i.e. that has a boundary but no holes)
can be mapped conformally to the 2-d unit disk
$\bD=\set{z \ \mid \ z=x+\ii y,\,\abs{z}\leq 1}$, with the boundary
of the disk corresponding to the boundary of $\cS$\begin{footnote}{We shall
restrict ourselves to this case
here, although our approach is more general; see \cite{LPD2010}).}\end{footnote}.
This mapping is called  ``conformally flattening''
\begin{footnote}{The uniformization theorem holds
for more general surfaces as well. For instance, surfaces without holes, handles or
boundaries can be mapped conformally to a sphere; if one point is
removed from such a surface, it can be mapped conformally to the
full plane. Surfaces with holes or handles can still be conformally
flattened to a piece of the plane.}\end{footnote}.
This flattening
process is accompanied by area distortion; the \emph{conformal
factor} $f(x,y)$ on the disk, varying from point to point,
indicates the area distortion factor produced by the operation.

One important practical implication of this theorem is that the
family of conformal maps between two surfaces can be characterized
naturally via the flattened representations of the surfaces: if
$\gamma$ is a conformal mapping from $\cS$ to $\cS'$, and $\varphi$
($\varphi'$) is a flattening (i.e. a conformal map to the disk
$\bD$) of $\cS$ ($\cS'$), then the family of all possible conformal
mappings from $\cS$ to $\cS'$ is given by $\gamma =  {\varphi'}
^{-1} \circ m \circ \varphi$, where $m$ ranges over all
the conformal bijective self-mappings of the unit disk $\bD$. We shall
call such $m$ \emph{disk-preserving M\"{o}bius
transformations}; they constitute a group,
the \emph{disk-preserving M\"{o}bius
transformation group} $\mathcal{M}$. Each $m$ in $\mathcal{M}$ is characterized by
3 parameters and given by the closed-form formula $m(z)=e^{\ii
\theta}(z-\alpha)(1-z\bar{\alpha})^{-1}$, where $\theta\in[0,2\pi)$,
$|\alpha|<1$. For our applications, it
is important that
the flattening process (starting from
a triangulated digitized version of $\cS$) and more importantly the
disk-preserving M\"{o}bius transformations can be computed fast and
with high accuracy; for more details, see
\cite{LF2009,LPD2010}.$\!$\begin{footnote} {If the digitization of
the surface is given as a point cloud, standard fast algorithms can
be used to determine an appropriate (e.g. Delauney) triangulation.
}\end{footnote} Note that the flattening map of a
surface $\cS$ is not unique; one can choose any arbitrary point of
$\cS$ to be mapped to the origin of the disk $\bD$, and any
direction through this point to become the ``$x$-axis''. The
transition from choosing one (center, direction) pair to another is
simply a disk-preserving M\"{o}bius transformation.
It is convenient to
equip the disk $\bD$ with its \emph{hyperbolic measure}
$d\eta(x,y)\,=\, [1-(x^2+y^2)]^{-2}dx\,dy$, invariant with respect
to M\"{o}bius transformations;
correspondingly, we set
$\Bf(x,y)= [1-(x^2+y^2)]^{2}\,f(x,y)$, so that
$\Bf\,d\eta\,=\,f\,dx\,dy\,$.
\subsection{Optimal mass transportation}
$\,$An integrable function $\mu$ is a (normalized) \emph{mass
distribution} on a domain $D$ if $\mu(u) \geq 0$ is well defined for
each $u \in D$, and $\int_D \mu(u)\,du \,=\,1$. If $\tau$ is a
differentiable bijection from $D$ to itself, the mass distribution
$\mu'=\tau_* \mu$ on $D$ defined by
$\mu(u)=\mu'(\tau(u))\,J_{\tau}(u)$ (where $J_{\tau}$ is the
Jacobian of the map $\tau$), is the \emph{transportation (or
push-forward) of $\mu$ by $\tau$} in the sense that, for any
arbitrary (non pathological) function $F$ on $D$, $\int_D F(u) \,
\mu'(u)\,du\,=\, \int_D F(\tau(u))\,\mu(u)\,du$. The total
\emph{transportation effort} is given by
$\mathcal{E}_{\tau}\,=\,\int_D d(u,\tau(u))\,\mu(u)\,du$,
where $d(u,v)$ denotes the distance between two points $u$ and $v$ in $D$. \\
If two mass distributions $\mu$ and $\nu$  on $D$ are given, then
the \emph{optimal mass transportation distance between $\mu$ and
$\nu$} (in the sense of Monge, see \cite{Villani2003}, p. 4) is the
infimum of the transportation effort $\mathcal{E}_{\tau}$, taken
over all the measurable bijections $\tau$
from $D$ to $D$ for which $\nu$ equals the transportation of $\mu$
by $\tau$. This set of bijections is hard to search; the
determination of an optimal mass transportation scheme becomes more
tractable if the mass ``at $u$'' need not all end up at the same end
point. One then considers measures $\pi$ on $D \times D$ with
marginals $\mu$ and $\nu$ (this means that for all continuous
functions $F,G$ on $D$,  $\int_{D\times D}F(u)\,d\pi(u,v)$
$=$ $\int_D F(u)\,\mu(u)\,du$ and $\int_{D\times D}G(v)\,d\pi(u,v)$
$=$ $\int_D G(v)\,\nu(v)\,dv$); the optimal mass transportation in this more general
Kantorovitch formulation is the infimum over all such measures $\pi$
of $E_{\pi}$ $=$ $\int_{D \times D}d(u,v)\,d\pi(u,v)$. A
comprehensive treatment
of optimal mass transport is in \cite{Villani2003}.
\section{2. New distances between 2-dimensional surfaces}
\subsection{Conformal Wasserstein distance (cW)}$\,$
One can use optimal mass transport to compare conformal factors
$\Bf$ and $\Bf'$ obtained by conformally flattening two surfaces,
$\cS$  and $\cS'$. If $m$ is a disk-preserving
M\"{o}bius transformation, then $\Bf$ and $m_*\Bf = \Bf \circ
m^{-1}$ are both equally valid conformal factors for $\cS$.
A standard approach to take
this into account is to ``quotient'' over
$\cM$, which leads to the
\emph{conformal Wasserstein distance}:
\begin{equation}
\label{eq:conf_Wass_dist}
\mbox{\hausaD}_{\mbox{\small{cW}}}(\cS,\cS')\!=\!\!\inf_{m \in
\cM}\!\left[\inf_{\pi \in \Pi(m_*
\Bf,{\Bf}^{'})}\!\int_{\bD\times\bD} \!\!\!\!\!\!\wt{d}\big
(z,z'\big ) d\pi (z,z')\right],
\end{equation}
where $\wt{d}(\cdot,\cdot)$ is the (conformally invariant)
hyperbolic distance\begin{footnote}{This is the geodesic distance
on $\bD$ induced by the hyperbolic Riemann metric tensor $d \eta$ on $\bD$. The geodesic
from the origin to any point $z$ in $\bD$ is the straight line
connecting them, and $\wt{d}(0,z)=\ln[(1+|z|)/(1-|z|)]$.}\end{footnote}
in $\bD$;
$\mbox{\hausaD}_{\mbox{\small{cW}}}$ satisfies then all the
properties of a metric \cite{LipDau2010}.
In particular, $\mbox{\hausaD}_{\mbox{\small{cW}}}(\cS,\cS')=0$ iff $\cS$
and $\cS'$ are isometric. However, computing this metric requires solving a
Kantorovitch mass-transportation problem for every candidate
$m$; even though the whole procedure has
polynomial runtime
complexity,
it is too heavy to be used
in practice for large datasets.
%
\subsection{Conformal Wasserstein neighborhood dissimilarity distance (cWn)}$\,$
We propose another natural way
to use Kantorovich's optimal mass transport
to compare surfaces $\cS$ and $\cS'$. Instead of determining
the most
efficient way to transport ``mass''
$\Bf$ from $z$ to $z'$, we can quantify how \emph{dissimilar} the
``landscapes'' are, defined by $\Bf$ and $\Bf'$ near $z$, resp.
$z'$, and replace the distance
$\wt{d}(\cdot,\cdot)$ by a \emph{measure of neighborhood
dissimilarity}. The neighborhood $N(0,R)$ around $0$ is given
by $N(0,R)=\{\,z\,;\,|z|<R \, \}$; neighborhoods around other points
are obtained by letting the disk-preserving M\"{o}bius
transformations act on $N(0,R)$:  for any $m$ in
$\mathcal{M}$ such that $z=m(0)$,  $N(z,R)$ is the image of $N(0,R)$ under the
mapping $m$. Next
we define the dissimilarity between $\Bf$ at $z$ and $\Bf'$ at $z'$:
\[
\bd^R _{\Bf,{\Bf}^{'}}(z,z')=\!\!\inf_{m \in \mathcal{M}, m(z)=z'}\!
\left[\,\int_{N(z,R)}|\Bf(w)-\Bf'(m(w))\,|\,d\eta(w)\right].
\]
It is straightforward to check that for all $m,m'$ in $\cM$, $\bd^R
_{m_* \Bf,{{m}^{'}}_* {\Bf}^{'} }(m(z),m'(z')) =\bd^R
_{\Bf,{\Bf}^{'}}(z,z')$. We now use optimal transport, and define
the \emph{conformal Wasserstein neighborhood dissimilarity distance}
between $\Bf$ and $\Bf'$:
\begin{equation}
\label{eq:conf_Wass_diss_dist}
\mathcal{D}^R_{\mbox{\small{cWn}}}(\cS,\cS')\,=\,\inf_{\pi \in
\Pi(\Bf,{\Bf}^{'})}\int_{\bD\yl{\times}
\bD}\bd_{\Bf,{\Bf}^{'}}(z,z')\,d\pi(z,z')~,
\end{equation}
where the superscript recalls that this definition depends on the
choice of the parameter $R$. For a proof that
this defines a true distance between (generic)
surfaces $\cS$ and $\cS'$, and further
mathematical properties, see
\cite{LipDau2010,LipDau2010tecrep}. One practical difference with
$\mbox{\hausaD}_{\mbox{\small{cW}}}$
is that (\ref{eq:conf_Wass_diss_dist})
requires solving only one Kantorovitch
mass-transportation problem once the special dissimilarity cost is
computed, resulting in a simpler optimization problem. To implement
the computation of these distances, we discretize the integrals and
the optimization searches, picking collections of discrete points on
the surfaces; the minimizing measure $\pi$ in the definition of
$\mathcal{D}^R_{\mbox{\small{cWn}}}(\cS,\cS')$ can then be used to
define a correspondence between points of $\cS$ and $\cS'$.
\subsection{Continuous Procrustes distance between surfaces (cP)}
$\,$ Both cW and cWn are \emph{intrinsic}: they use only
information ``visible'' from within each surface, such as geodesic
distances between pairs of points; consequently they do not distinguish
a surface from any of its isometric embeddings in 3-D.
The continuous Procrustes (cP) distance \cite{LRFD2010procrustes}
described in this section
uses some {\em extrinsic} information as well; it fails to distinguish two
surfaces only if one is obtained by applying to the other a {\em rigid} motion
(which is a very special isometry).

The {\em (standard) Procrustes distance} is defined between discrete
sets of points $\bX=(X_n)_{n=1,\ldots,N} \subset \cS$ and
$\bY=(Y_n)_{n=1,\ldots,N} \subset \cS'$ by
minimizing over all rigid motions:
$d_P(\bX,\bY)\,=\,\min_{R\mbox{\small{ rig.mot.}}\!} \left[\,
\left(\,\sum_{n=1}^N|R(X_n)-Y_n|^2\,\right)^{1/2}\,\right]$, where
$\abs{\cdot}$ denotes the standard Euclidean norm.
\footnote{It is interesting to note that in \cite{Mem_again}, a Kantorovich version of $d_P$ was
introduced, and
its equivalence to the Gromov-Wasserstein distance (when the shapes are endowed with Euclidean distances)
was proved.}
Often
$\bX,\, \bY$ are sets of {\em landmarks} on two surfaces, and $d_P(\bX,\bY)$
is interpreted as a distance between these
surfaces.
This practice has several drawbacks: 1) $d_P(\bX,\bY)$
depends on the (subjective) choice of $\bX,\bY$, which makes it a
not necessarily ``well-defined'' or easily reproducible proxy for a surface distance;
2) the (relatively) small number of $N$ landmarks on each surface disregards
a wealth of geometric data; 3) identifying and recording the $x_n,\,y_n$ is
time consuming and requires expertise.

We eliminate all these drawbacks by a
{\bf landmark-free} approach, introducing the
{\em continuous Procrustes distance}. Instead
of relying on experts to identify ``corresponding'' discrete
subsets of $\cS$ and $\cS'$,
we consider a family of
\emph{continuous maps} $a:\cS\rightarrow \cS'$ between the surfaces,
and rely on optimization to identify
the ``best'' $a$. The earlier exact correspondence of
one point $Y_n$ to one point $X_n$, and the
(tacit) assumption that $\bX$ ($\bY$) collectively represent all the noteworthy
aspects of $\cS$ ($\cS'$) in a balanced way,
are recast as requiring that the ``correspondence map''
be \emph{area-preserving} \cite{LRFD2010procrustes},
that is, for every (measurable)
subset $\Omega$ of $\cS$,
$\int_{\Omega}dA_{\cS}=\int_{a(\Omega)}dA_{\cS'}$, where $dA_{\cS}$
and $dA_{\cS'}$ are the area elements on the surfaces induced by their embeddings in $\R^3$.
We denote
$\cA(\cS,\cS')$ the set of all these area-preserving diffeomorphisms.
For each $a$ in
$\cA(\cS,\cS')$, we set $\bd(\cS,\cS',a)^2= \min_{R\mbox{\small{
rig.mot.}}\!}\int_{\cS}|R(x)-a(x)|^2dA_{\cS}$; the \emph{continuous
Procrustes distance between $\cS$ and $\cS'$} is then
\begin{equation}
\label{eq:cont_Procr_dist} \bD_P(\cS,\cS')=\inf_{a \in
\cA(\cS,\cS')} \bd(\cS,\cS',a)~.
\end{equation}
This defines a metric distance on the
space of surfaces (up to rigid motions: for congruent surfaces the
distance is 0) \cite{LRFD2010procrustes}. Minimizing over rigid
motions is easy; there exist closed form formulas, as in
the discrete case. But the second set over which to minimize,
$\mathcal{A}(\cS,\cS')$, is an unwieldy, formally
infinite-dimensional manifold, hard to explore\begin{footnote}{It
can be viewed as the continuous analog to
the exponentially large
group of permutations.}\end{footnote}. This is an
optimal transport problem again, now in the much harder Monge
formulation. For ``reasonable'' surfaces (e.g., surfaces
with uniformly bounded curvature), transformations $a$ close to
optimal are close to conformal
\cite{LRFD2010procrustes}. This crucial insight allows limiting the
search to the much smaller space of maps
obtained by small deformations
of conformal
maps. Concretely, we
compose a conformal map (represented as a
M\"{o}bius transformation) $m\in\mathcal{M}$  with maps
$\chi$ and $\varrho$, where $\varrho$ is a smooth map
that roughly aligns high density peaks, and $\chi$ is a special
deformation (following
\cite{Moser}) using local diffusion to make
$\chi\circ \varrho \circ m$
area preserving (up to approximation error).  For each choice of peaks $p,\,p'$
in the
conformal factors of $\cS$, $\cS'$, the algorithm 1) runs through
the 1-parameter family of $m$ that map $p$ to $p'$;
2) constructs a
map $\varrho$ that aligns the other peaks, as best possible;
and 3) computes $\bd(\cS,\cS',\varrho\circ m)$.
Repeat for all choices of $p,\,p'$; the $\varrho\circ m$ that minimizes $\bd$ is then
deformed to be area preserving, producing
the map $a = \chi\circ \varrho \circ m$;
$\bd(\cS,\cS',a)$ and $a$ are our approximate
$\bD_P(\cS,\cS')$ and correspondence map, respectively.
(More in Supplementary Materials.)

\vspace*{-.4 cm}

\section{3. Application to anatomical datasets}
To test
our approach, we used three
independent datasets, representing three different regions of the
skeletal anatomy, of humans, other primates, and their close
relatives. Digitized surfaces
were obtained
from High Resolution X-ray Computed Tomography ($\mu$CT) scans (see Supplementary Materials) of $\,$(A) 116
second mandibular molars of prosimian primates and non-primate close
relatives, and $\,$(B) 57
proximal first metatarsals of prosimian primates, New
and Old World monkeys, and $\,$(C) 45 distal radii of apes and humans.
For every pair
of surfaces, the output of our algorithms
consists of 1) a correspondence map for the {\em whole} surface
(i.e. not
just a few points), and 2) a non-negative number
giving their dissimilarity (where zero means they are
\emph{isometric} or \emph{congruent}).
Typical running times for a pair of surfaces were $\sim 20$
sec. for cP, $\sim 5$ min. for cWn.
To evaluate the performance
of the algorithms, we compared the
outcomes to those determined independently by morphologists.
Using the same set of digitized surfaces, geometric
morphometricians collected landmarks on each, in the conventional
fashion \cite{Mitteroeker_Gunz2009}, choosing them to reflect
correspondences considered biologically and evolutionarily
meaningful (see Supplementary Materials). These landmarks
determine ``discrete'' Procrustes distances for every two
surfaces (see $\mathsection$ 2), here called Observer-Determined Landmarks Procrustes
(ODLP) distances. For each of the
three distances we obtain thus a (symmetric) matrix. \\
\vspace*{-.6cm}

\subsection{Comparing the distance (dissimilarity) matrices}
$\,$ We compare cWn- and cP- with ODLP-matrices in two different ways.
Sets of
distances are far from independent,
and it is traditional to assess
the relationship between distance matrices by a Mantel correlation
analysis \cite{Mantel1967}:
first correlate the entries in the two square arrays, and then compute the
fraction, among all possible relabelings of the rows/columns for one of them,
that leads to a larger correlation coefficient; this {\em Mantel significance} is
a stronger indicator than the correlation coefficient itself.
Table 1 gives the results for our datasets.
\begin{figure}[b]
$~$
\vspace*{-0.9 cm}

\begin{center}
\includegraphics[height=2.2 cm,width=1.02 \columnwidth]{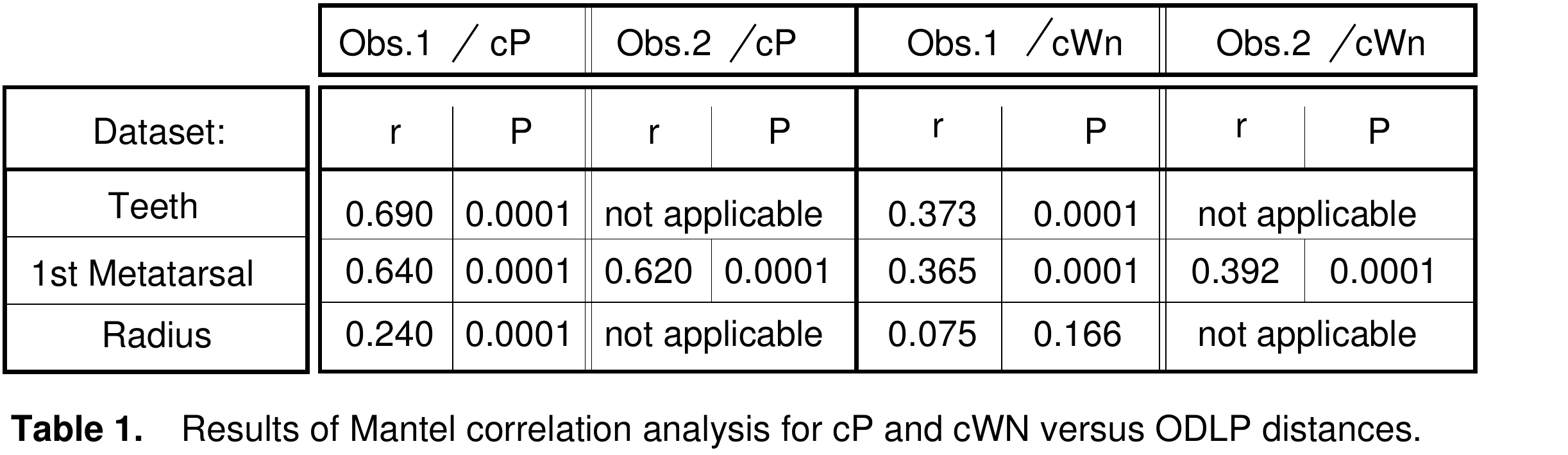}
\end{center}
$~$
\vspace*{-.4 cm}

\end{figure}

\noindent
In all cases the Mantel significance between ODLP and cP distances
is higher than that between ODLP and cWn. This

indicates that
distances computed using cP match those determined by
morphometricans' better than those using cWn.

Figure 1 illustrates the relationship of cP, cWn, and ODLP distances
in a different way.
In each of the
two square matrices (corresponding to cP and cWn, each vs. ODLP),
the color of each pixel indicates the value of the entry (using a red-blue colormap,
with deep blue representing 0, and saturated red the largest value);
upper right triangular halves correspond to cP or cWn;
(identical) lower left halves to ODLP. The same ordering of samples
is used in the three cases, with samples ordered so that nearby
samples typically have smaller distances. This type of display is
especially good to compare the structure of two distance matrices
for small distances, often the most reliable
\begin{footnote}{In some modern data analysis methods, such as
diffusion-based or graph-Laplacian based methods, only the small
distances are retained, to be used in spectral methods that ``knit''
the larger scale distances of the dataset together more
reliably.}\end{footnote}.
Note the better symmetry along the diagonal for ODLP/cP comparison
on the left: in this comparison, as in the previous one, cP outperforms cWn.
\vspace*{-0.2cm}

\begin{figure}[h]
\centering 
\includegraphics[width=\columnwidth]{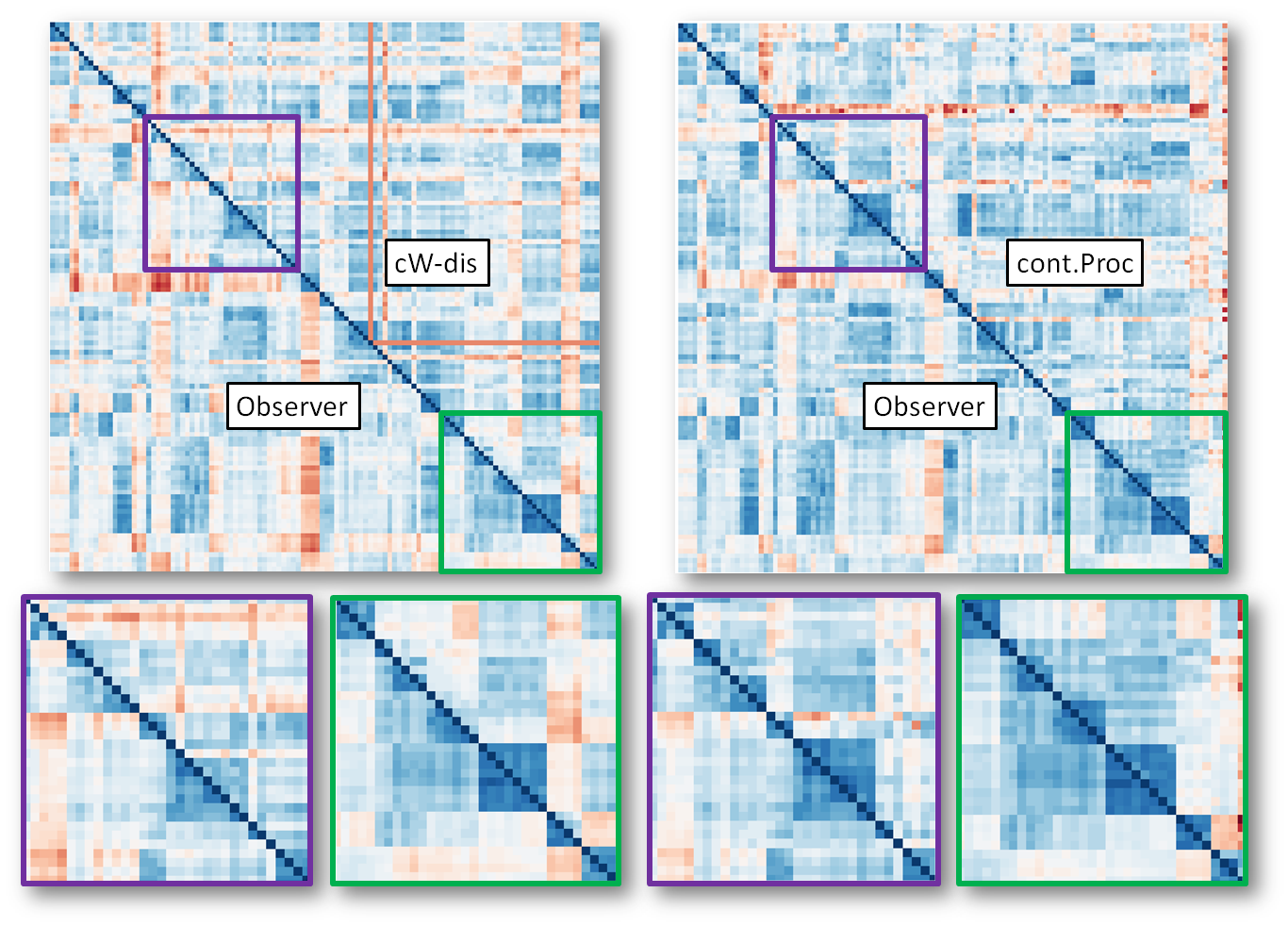}
\caption{For small distances,
the structures of the matrices with cP-, cWn-distances and distances
based on Observer Landmarks (ODLP) are very similar, with cP  (on the
right) the most similar to ODLP. The dataset illustrated here is
dataset (A).}
\label{fig:matrices}
\vspace*{-0.3cm}

\end{figure}%
\vspace*{-0.9cm}

\subsection{Comparing scores in taxonomic classification}
$\,$Accurately placed ODL usually result in smaller ODLP distances
between specimens representing individuals of the same species/genus
than between individuals of different species/genera.
\begin{figure}[htb]
$~$
\vspace*{-0.7 cm}

\begin{center}
\includegraphics[height=2.2 cm,width=1 \columnwidth]{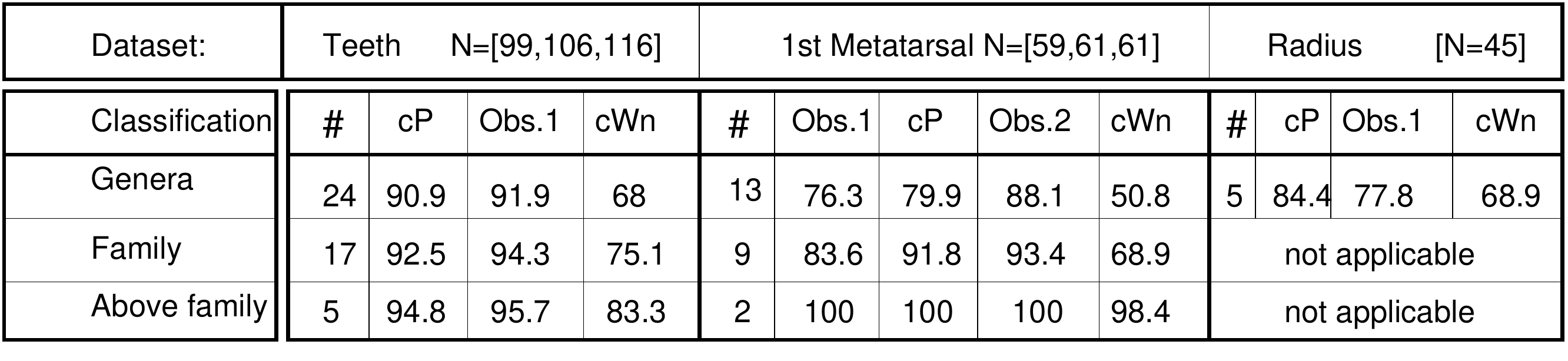}
\end{center}
$~$
\vspace*{-.8 cm}

\end{figure}
\noindent
To assess whether this holds as well for the algorithmic cP and cWn
distances, we run three taxonomic
classification analyses
on each data set, one using ODLP distances, and two using
cP and cWn distances\begin{footnote}
{We do {\bf not}
claim this is a new/alternative method for automatic species or genus
identification. Reliable automatic species recognition uses, in
addition, auditory, chemical, non-geometric morphological and other
data, analyzed by a range of methods; see e.g.
\cite{WEEKS_ET_AL1999,MacLeod2007, MacLeod2008} and references.
Comparison of taxonomic classification based on
human-expert-generated vs. algorithm-computed-distances is
meant only as a quantitative evaluation based on biology rather than
mathematics. }\end{footnote}, with
a ``leave one
out''  procedure: each specimen (treated as
unknown) is assigned to the taxonomic group of its nearest neighbor
among the remainder of the specimens in the data set (treated as known).
Table 2 lists success rates (in \%) for three different classification
queries for the three datasets. For each dataset $N$ is
the number of objects; for each query $\#$ is the number of groups.
Classifications based on the cP distances are
similar in accuracy to those based on the ODLP distances,
outperforming the cWn distances for all three of our anatomic datasets.

\noindent
{\bf Note:} a similar classification based on topographic variables is less accurate:
for the 99 teeth belonging to 24 genera, only 54 (54 \%) were classified correctly
with a classification based on the four topographic variables
Energy, Shearing Quotient, Relief Index, OPC.
(Details in Supplementary Materials.)
\vspace*{-.5 cm}

\subsection{Comparing the correspondence maps}
$\,$ Morphometric analyses are based on the identification of
{\em corresponding landmark points} on each of $\cS$ and $\cS'$; the cP
algorithm constructs a {\em correspondence map} $a$ from $\cS$ to $\cS'$.
(The correspondence induced by cWn is less smooth and will not be considered here.)
For each landmark point $L$ on
$\cS$, we can compare the location on $\cS'$
of its images $a(L)$
with the location of the corresponding landmark points $L'$.
Fig. 2 shows that the ``propagated''  landmarks
$a(L)$ typically turn out to be very close to those of the ``true'' landmarks
$L'$.
(More in Supplementary Materials.)
\begin{figure}[h]
\centering
\includegraphics[width=\columnwidth]{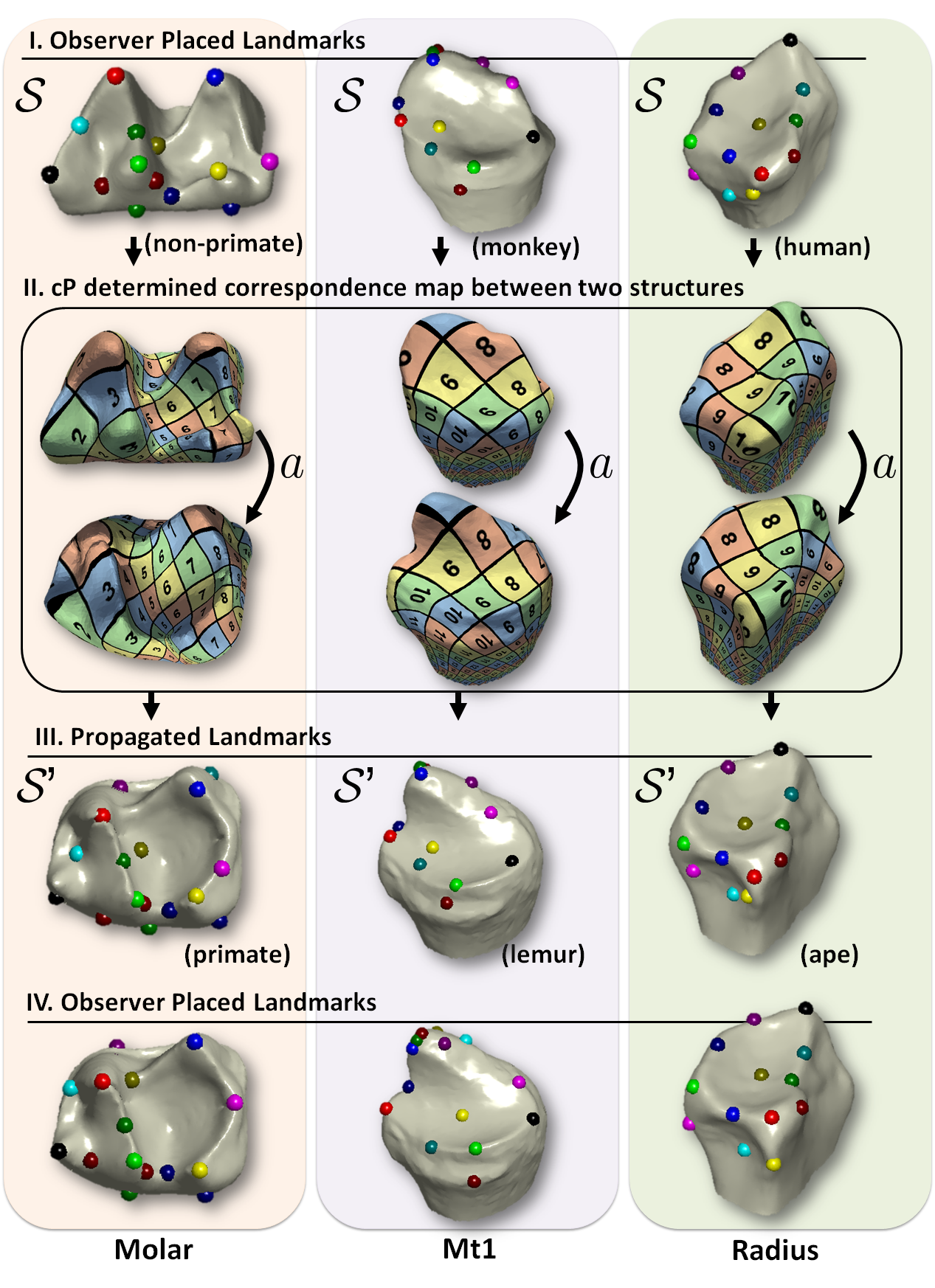}
\caption{Observer-placed landmarks can be propagated from structure
M (A) using cP determined correspondence maps (B) to another specimen
N (C). The similarity between propagated landmarks in (C) and
observer placed landmarks in (D) on N shows the success of the
method, and makes explicit the geometric basis for the
observer-determinations.}
\label{fig:correspondences}
\vspace*{-0.2cm}

\end{figure}
%
\subsection{An application}
These comparisons show our algorithms capture
biologically informative shape
variation.
But scientists are
interested in more than overall shape!
We illustrate how correspondence maps could be
used to analyze more specific features. In comparative
morphological and phylogenetic studies,
anatomical identification of certain features (e.g., particular
cusps on teeth) is controversial in some cases; an example of this is
the distolingual corner of
sportive lemur (\emph{Lepilemur}) lower molars in our data
set (A)
\cite{Swindler2002,Schwartz_Tat1985}.

In such controversial cases, transformational homology
\cite{Patterson1982} hypotheses are usually supported by a specific
comparative sample or inferred morphocline
\cite{Osborn1907,VanValen1982,Schwartz_Tat1985}. \emph{Lepilemur} is
thought by some researchers to lack a cusp known as an entoconid
(Fig. 3) but to have a hypertrophied metastylid cusp that ``takes
the place'' of the entoconid \cite{Schwartz_Tat1985} in other taxa.
Yet, in comparing a \emph{Lepilemur} tooth to a more ``standard''
primate tooth, like that of \emph{Microcebus}, both seem to have the
same basic cusps; alternatives to the viewpoint of
\cite{Schwartz_Tat1985} have therefore also been argued in the
literature \cite{Swindler2002}. However, another lemur,
\emph{Megaladapis} (now extinct), arguably a closer relative of
\emph{Lepilemur} than \emph{Microcebus}, has an entoconid that is
very small and a metastylid that is rather large, thus providing an
evolutionary argument supporting the original hypothesis. (For more
details, see Supplementary Materials.)
%
Such arguments can now be made more precise.
We can propagate (as in Fig. 2) landmarks from the
\emph{Microcebus} to the \emph{Lepilemur} molar; this direct
propagation matches the entoconid cusp of \emph{Microcebus} with the
controversial cusp of \emph{Lepilemur} (Fig. 3, path 1), supporting
\cite{Swindler2002}.
In contrast, when we propagate
landmarks in different steps, either from \emph{Microcebus} to
\emph{Megaladapis} and then to \emph{Lepilemur}(Fig. 3, path
2), or through the extinct \emph{Adapis} and extant \emph{Lemur}
(Fig. 3, path 3), the
\emph{Lepilemur} metastylid takes the place of the \emph{Microcebus}
entoconid, supporting \cite{Schwartz_Tat1985}.
Automatic propagation of landmarks via mathematical algorithms
recenters the controversy on the (different)
discussion of which propagation channel is most suitable.
\begin{figure}[h]
\centering
\includegraphics[width=\columnwidth]{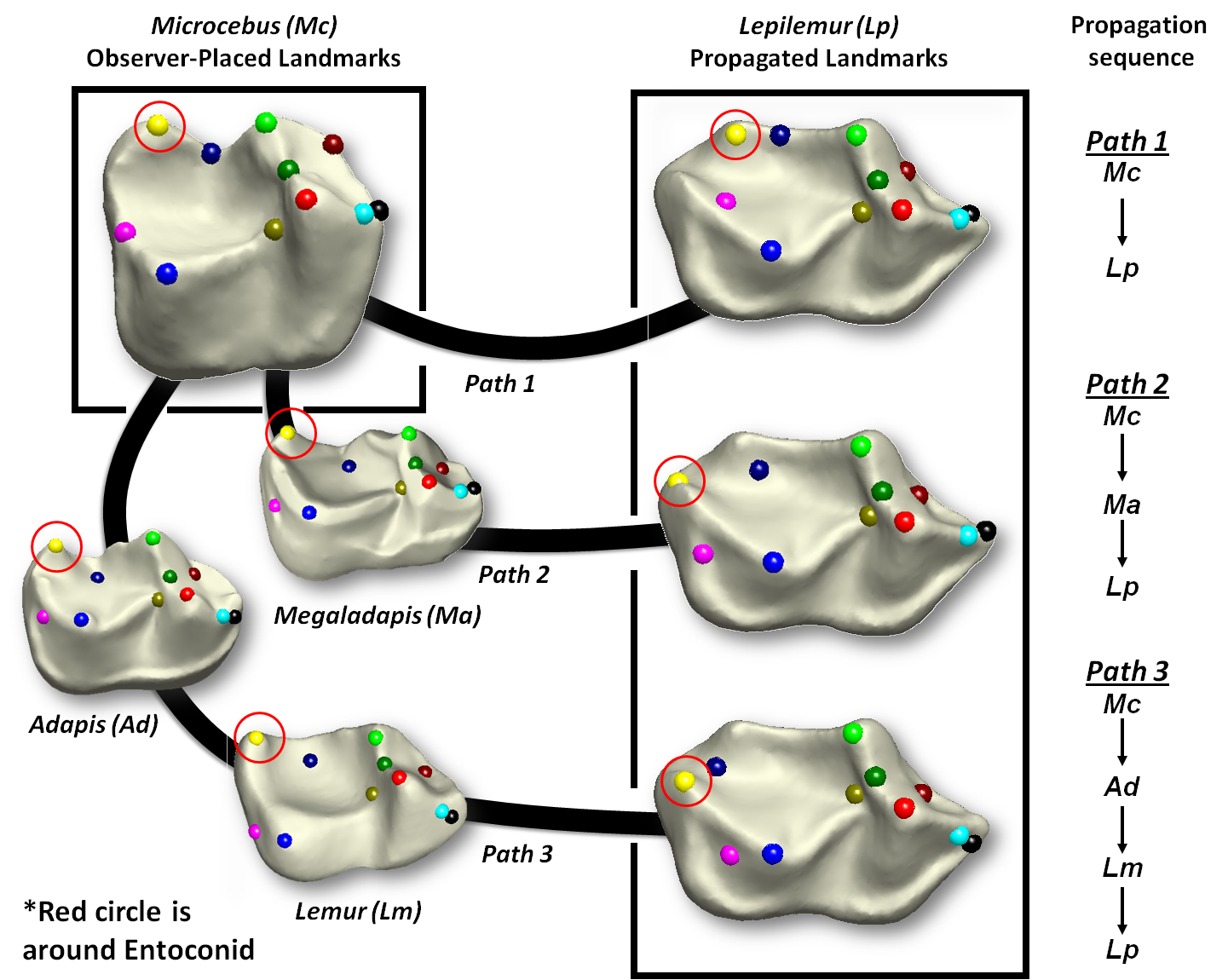}
\caption{Observer-placed landmarks on a tooth of \emph{Microcebus}
are propagated using cP-determined correspondence maps to a tooth of
\emph{Lepilemur}. Path 1 is direct, Path 2 and 3 have intermediate
steps, representing step-wise propagation between teeth of other
taxa.}
\label{fig:pathways}
\end{figure}
%
\subsection{Summary and Conclusion}
$\,$
New distances between 2-D surfaces, 
with fast numerical implementations, were shown to lead to
fast,
landmark-free algorithms that map anatomical
surfaces automatically to other instances of anatomically equivalent
surfaces, in a way that mimics
accurately the detailed
feature-point correspondences recognized qualitatively by
scientists, and that preserves information on taxonomic structure
as well as observer-determined-landmark distances.
Moreover, the correspondence maps
thus generated can incorporate, in their tracking of point features,
evolutionary relationships inferred to link different taxa together.

Our approach makes
morphology
accessible to non-specialists and allows the documentation
of anatomical variation and quantitative traits with previously
unmatched comprehensiveness and
objectivity. More
frequent, rapid, objective, and comprehensive construction of
morphological datasets will allow the study of morphological
diversity's evolutionary significance to be better synchronized with
studies incorporating genetic and developmental information, leading
to a better understanding of anatomical form and its genetic basis,
as well as
the evolutionary processes that have contributed to their diversity
among living things on earth.


\begin{acknowledgments}
Access was provided by 1) curators and staff at the Am. Museum of
Natl. Hist. (to dental specimens to be molded, cast, and
microCT-scanned), 2) D. Krause, J. Groenke of Stony Brook Univ.
Dept. of Anatomical Sciences (to facilities of the Vertebrate
Paleontology Fossil Preparation Lab.), and
3) C. Rubin, S. Judex of Stony Brook Univ. Dept. of Biomedical Engineering and the Center for Biotechnology (to microCT scanners, for digital imaging of tooth casts). \\
An NSF DDIG through Physical Anthropology, the Evolving Earth
Foundation and the American Society of Physical Anthropologists
provided funding to DMB; YL, JP, TF, and ID were supported by NSF
and AFOSR grants.
\end{acknowledgments}





\end{article}








\end{document}